\documentclass{elsart}

\usepackage{amsmath,amssymb,amsbsy,amsfonts,  latexsym,amsopn,amstext, amsxtra,euscript,amscd}

\newcommand{\set}[1]{{\left\{#1\right\}}}    
\newcommand{\K}{{K}}
\def\C{\mathbb C}
\def\Q{\mathbb Q}
\def\F{\mathbb F}

\def\Z{\mathbb Z}
\def\O{\mathcal O}
\def\R{\mathcal R}

\def\La{\Lambda}
\def\Th{\theta}
\def\th{\theta}

\def\om{\omega}

\def\<{<}
\def\>{>}
\def\O{\mathcal O}

\def\L{\Lambda}

\def\proof{\textbf{Proof: }}

\begin{document}

\begin{frontmatter}



\title{Codes over rings of size $p^2$ and lattices over imaginary quadratic fields}


\thanks[label1]{Partially supported by a NATO grant}

\author{T. Shaska\thanksref{label1}}
\address{Department of Mathematics and Statistics, \\ Oakland University,  Rochester, MI, 48309. USA
\\ and \\
Department of Computer Science and Electrical Engineering, \\
University of Vlora, \\
Vlora, ALBANIA } \ead{shaska@oakland.edu}

\author{C. Shor}
\address{Department of Mathematics, \\ Bates College,  3 Andrews Road, Lewiston, ME, 04240.}
\ead{cshor@bates.edu}

\author{S. Wijesiri}
\address{Department of Mathematics and Statistics, \\ Oakland University,  Rochester, MI, 48309.}
\ead{gswijesi@oakland.edu}

\begin{abstract}
Let $\ell>0$ be a square-free integer congruent to 3 mod 4 and $\O_K$ the ring of integers of the
imaginary quadratic field $K=Q(\sqrt{-\ell })$. Codes $C$ over rings $\O_K / p \O_K$ determine lattices
$\Lambda_\ell (C) $ over $K$. If $ p \nmid \ell$ then the ring $\R:=\O_K / p \O_K$ is isomorphic to
$\F_{p^2}$ or $\F_p \times \F_p$. Given a code $C$ over $\R$, theta functions on the corresponding
lattices are defined. These theta series $\theta_{\Lambda_{\ell}(C)}$  can be written in terms of the
complete weight enumerator of $C$. We show that for any two $\ell < \ell^\prime$ the first $\frac {\ell +
1} 4$ terms of their corresponding theta functions are the same. Moreover, we conjecture that for $\ell >
\frac {p(n+1)(n+2)} 2$ there is a unique complete weight enumerator corresponding to a given theta
function. We verify the conjecture for primes $p< 7$ and $\ell \leq 59$.
\end{abstract}


\begin{keyword}
codes \sep  lattices \sep  theta functions
\end{keyword}

\end{frontmatter}

\section{Introduction}
Let $\ell>0$ be a square-free integer congruent to 3 modulo 4,  $K=\Q(\sqrt{-\ell })$ be the imaginary
quadratic field, and $\O_K$ its ring of integers. Codes, Hermitian lattices, and their theta-functions over
rings $\R:=\O_K / p \O_K$, for small primes $p$, have been studied by many authors, see \cite{ba, ms1,
ms2},  among others. In \cite{ba}, explicit descriptions of theta functions and MacWilliams identities are
given for $p=2, 3$. In \cite{SS} we explored codes $C$ defined over $\R$ for $p >2$. For any $\ell$ one
can construct a lattice $\Lambda_\ell (C)$ via Construction A  and define theta functions based on the
structure of the ring $\R$. Such constructions suggested some relations between the complete weight
enumerator of the code and the theta function of the corresponding lattice. In this paper we give complete
proofs of some of the theorems in \cite{SS}.   Furthermore, we study the weight enumerators of such
codes in terms of the theta functions of the corresponding lattices.  This paper is organized as follows.

In section 2 we give a brief overview of the basic definitions for codes and lattices and define theta
functions over $\F_p$. We define the theta series $\th_{\Lambda_{a, b} } (q)$ for all cosets in $p\O_K$ and
determine relations among such theta series. Two such theta series $\th_{\Lambda_{a, b} } (q)$ and
$\th_{\Lambda_{m, n} } (q)$ are the same when $(m, n)$ is congruent modulo $p$ to one of the ordered pairs
$(a, b), (-a-b, b), (-a, -b), (a+b, -b)$. This implies that we have at most $\frac {(p+1)^2} 4$ theta series,
and when $\ell > 12p^2+1$ we have exactly $\frac {(p+1)^2} 4$ theta series. In section 3 we define
theta functions on the lattice defined over $\R:=\O_K / p \O_K$.  We prove in \cite{SS} that such a theta
series is equal to the evaluation of the complete weight enumerator of the code on the theta series of cosets of $p\O_K$.

In section 4, we address a special case of a general problem of the construction of lattices: the injectivity of
Construction A.  For codes defined over an alphabet of size four (regarded as a quotient of the ring of
integers of an imaginary quadratic field), the problem is solved completely in \cite{SV}.  The analogous
questions are asked for codes defined over  $\F_{p^2}$ or $\F_p \times \F_p$. The main obstacle seems
to be expressing the theta function in terms of the symmetric weight enumerator of the code. However, the
theta function  $\theta_{\Lambda_{\ell}(C)}$ can be expressed in terms of the complete weight
enumerator of the code. We expect that similar results as for $p=2$ hold also for odd primes. However,
we are not able to get explicit bound for $p>2$. In section 5 we display some computational results for
$p=3$. Such results confirm our results of section 4. We compute the theta series for $p=3$, $n=3, 4, 5$,
and $\ell \leq 59$. We conjecture that for $\ell > \frac {p (n+1)(n+2)} 2$ for each given theta series exists
at most one complete weight enumerator polynomial corresponding to this theta series.

\section{Preliminaries}
Let $\ell>0$ be a square free integer and $K=\Q(\sqrt{-\ell })$ be the imaginary quadratic field with
discriminant $d_K$. Recall that $ d_K=   -\ell$ if $\ell\equiv 3 \mod 4,$ and $d_K=  -4\ell$ otherwise. Let
$\O_K$ be the ring of integers of $K$. A  lattice $\La$ over $K$ is an $\O_K$-submodule of $K^n$ of full
rank.  The Hermitian dual is defined by
\begin{equation}
\L^* = \set{x \in \K^n \; | \; x \cdot \bar{y} \in \O_K, \text{for all } \, \,  y \in \La},
\end{equation}
where $ x \cdot y := \sum_{i=1}^n x_i y_i$ and $\bar{y}$ denotes component-wise complex conjugation. In the
case that $\L$ is a free $\O_K$ - module, for every $\O_K$ basis $ \set {v_1 , v_2 , ...., v_n}$ we can
associate a Gram matrix G($\L$) given by  $G(\L) = (v_i.v_j)_{i,j=1}^n $ and the determinant $\det \L
:=\det(G) $ defined up to squares of units in $\O_K$. If $ \L = \L^* $ then $ \L$ is Hermitian self-dual (or
unimodular) and integral if and only if $\L \subset \L^*$. An integral lattice has the property $   \L
\subset \L^* \subset \frac{1}{det \L} \L$. An integral lattice is called even if $x\cdot x\equiv 0 \mod 2$
for all $x\in\L$, and otherwise it is odd. An odd unimodular lattice is called a Type 1 lattice and even
unimodular lattice is called a Type 2 lattice.

The theta series of a lattice $\L$ in $K^n$ is given by
\[ \Th_\L(\tau) = \sum_{z \in \L} e^{\pi i\tau z\cdot\bar{z}},\]
where $ \tau \in H =\set { z \in \C : Im(z)>0}.$ Usually we let $q = e^{\pi i \tau}. $ Then, $ \Th_\L(q) =
\sum_{z \in \L} q^{z\cdot\bar{z}}$. The one dimensional theta series (or Jacobi's theta series) and its
shadow are given by
\[\th_3(q)=\sum_{n\in\Z} q^{n^2}, \quad \th_2(q)=\sum_{n\in\Z} q^{(n+1/2)^2}=\sum_{n\in\Z+\frac{1}{2}}q^{n^2}.\]
Let $\ell\equiv 3 \mod 4$ and $d$ be a positive number such that $\ell=4d-1$. Then, $-\ell\equiv 1 \mod 4$.
This implies that the ring of integers is $\O_K=\Z[\om_\ell]$, where $\om_\ell=\frac{-1+\sqrt{-\ell}}{2}$ and
$\om_\ell^2 + \om_\ell+d=0$. The principal norm form of $K$ is given by
\begin{equation}\label{eq1}
Q_d(x,y) = |x-y\om_\ell|^2 = x^2+xy+dy^2.
\end{equation}
The structure of $\O_K/p\O_K$ depends on the value of $\ell$ modulo $p$. For $(\frac{a}{p})$ the Legendre
symbol,
\begin{equation}
 \O_K/p\O_K =
 \begin{cases}
   \F_p\times\F_p & \text{if $(\frac{-\ell}{p})=1$}, \\
   \F_{p^2} & \text{if $(\frac{-\ell}{p})=-1$}, \\
   \F_p+u\F_p \text{ with $u^2=0$} & \text{if $p\mid\ell$}.
 \end{cases}
\end{equation}
In this paper we will focus on the cases when $p\nmid\ell$.
\subsection{Theta functions over $\F_p$}
Let $q=e^{\pi i\tau}$.  For integers $a$ and $b$ and a prime $p$, let $\La_{a,b}$ denote the coset
$a-b\om_\ell+p\O_K$.  The theta series associated to this coset is
\begin{equation}
\begin{split}
\th_{\La_{a,b}}(q)  &=  \sum_{m,n\in\Z}q^{|a+mp-(b+np)\om_\ell|^2} \\
& =  \sum_{m,n\in\Z}q^{Q_d(mp+a,np+b)}  =  \sum_{m,n\in\Z}q^{p^2Q_d(m+a/p, n+b/p)}.
\end{split}
\end{equation}
For a prime $p$ and an integer $j$, consider the one-dimensional theta series
\begin{equation}\label{eqn:one-dim-series}
\th_{p,j}(q):=\sum_{n\in\Z}q^{(n+j/2p)^2}.
\end{equation}
Note that $\th_{p,j}(q)=\th_{p,k}(q)$ if and only if $j\equiv\pm k \mod 2p$.
\begin{lem}\label{lemma:in-terms-of-one-dim-series}
One can write $\th_{\La_{a,b}}(q)$ in terms of one-dimensional theta series defined above in Equation \ref{eqn:one-dim-series}. In particular,
\begin{equation}
\th_{\La_{a,b}}(q)=\th_{p,b}(q^{p^2\ell})\th_{p,2a+b}(q^{p^2})+\th_{p,b+p}(q^{p^2\ell})\th_{p,2a+b+p}(q^{p^2}).
\end{equation}
\end{lem}

\proof We use the fact that $Q_d(m,n)=m^2+mn+dn^2=(m+\frac{n}{2})^2+\frac{\ell n^2}{4}$.
\begin{small}
\begin{equation*}
\begin{split}
\th_{\La_{a,b}}(q) & =  \sum_{m,n\in\Z}q^{Q_d(mp+a,np+b)} \\
& =  \sum_{m,n\in\Z}q^{\left(mp+a+\frac{np+b}{2}\right)^2+\frac{\ell(np+b)^2}{4}} \\
& =  \sum_{n\in\Z} q^{\ell\frac{(np+b)^2}{4}}\sum_{m\in\Z}q^{\left(mp+a+\frac{np+b}{2}\right)^2} \\
& =  \sum_{n\in\Z} q^{\ell p^2(\frac{n}{2}+\frac{b}{2p})^2}\sum_{m\in\Z}q^{p^2(m+\frac{2a}{2p}+\frac{n}{2}+\frac{b}{2p})^2} \\
& = \sum_{n \text{ even}}q^{\ell
p^2(\frac{n}{2}+\frac{b}{2p})^2}\sum_{m\in\Z}q^{p^2(m+\frac{n}{2}+\frac{2a+b}{2p})^2}+
\sum_{n \text{ odd}}q^{\ell p^2(\frac{n}{2}+\frac{b}{2p})^2}\sum_{m\in\Z}q^{p^2(m+\frac{n}{2}+\frac{2a+b}{2p})^2} \\
& =  \th_{p,b}(q^{p^2\ell})\th_{p,2a+b}(q^{p^2})+\th_{p,b+p}(q^{p^2\ell})\th_{p,2a+p+b}(q^{p^2}).
\end{split}
\end{equation*}
\end{small}
This completes the proof. \qed

It would be interesting to determine what happens to the distribution of points on these cosets as $\ell$
increases. In other words, is there any relation among $\th_{\La_{a,b}}(q)$ as $\ell $ increases?
\begin{lem}\label{lemma-congruences}
For any integers $a,b,m,n$, if the ordered pair $(m,n)$ is congruent modulo $p$ to one of $(a, b), (-a-b, b),
(-a, -b), (a+b, -b)$, then $\th_{\La_{m,n}}(q)=\th_{\La_{a,b}}(q)$.
\end{lem}

\proof We aim to find sufficient conditions on $a,b,m,n$ so that $\th_{\La_{m,n}}(q)=\th_{\La_{a,b}}(q)$.
By Lemma \ref{lemma:in-terms-of-one-dim-series},
\[\th_{\La_{m,n}}(q)=\th_{p,n}(q^{p^2\ell})\th_{p,2m+n}(q^{p^2})+\th_{p,n+p}(q^{p^2\ell})\th_{p,2m+n+p}(q^{p^2})\]
and
\[\th_{\La_{a,b}}(q)=\th_{p,b}(q^{p^2\ell})\th_{p,2a+b}(q^{p^2})+\th_{p,b+p}(q^{p^2\ell})\th_{p,2a+b+p}(q^{p^2}).\]
In particular, if we have
\begin{eqnarray}\label{eqn:case1-first-eqn}
& \th_{p,n}(q^{p^2\ell}) &= \th_{p,b}(q^{p^2\ell}) \\
& \th_{p,2m+n}(q^{p^2}) &= \th_{p,2a+b}(q^{p^2}), \\
& \th_{p,n+p}(q^{p^2\ell}) &= \th_{p,b+p}(q^{p^2\ell})\\
& \th_{p,2m+n+p}(q^{p^2}) &= \th_{p,2a+b+p}(q^{p^2}), \label{eqn:case1-last-eqn}
\end{eqnarray}
(that is, equating the first terms, equating the second terms, etc.) then we will have $\th_{\La_{m,n}}(q)=\th_{\La_{a,b}}(q)$.

Similarly, if we change the order of the terms in $\th_{\La_{a,b}}(q)$ to obtain
\[\th_{\La_{a,b}}(q)=\th_{p,b+p}(q^{p^2\ell})\th_{p,2a+b+p}(q^{p^2})+\th_{p,b}(q^{p^2\ell})\th_{p,2a+b}(q^{p^2}),\]
we will have $\th_{\La_{m,n}}(q)=\th_{\La_{a,b}}(q)$ if
\begin{eqnarray}\label{eqn:case2-first-eqn}
& \th_{p,n}(q^{p^2\ell}) & = \th_{p,b+p}(q^{p^2\ell}) \\
& \th_{p,2m+n}(q^{p^2}) & = \th_{p,2a+b+p}(q^{p^2}), \\
& \th_{p,n+p}(q^{p^2\ell}) & = \th_{p,b}(q^{p^2\ell}) \\
& \th_{p,2m+n+p}(q^{p^2}) & = \th_{p,2a+b}(q^{p^2}).\label{eqn:case2-last-eqn}
\end{eqnarray}

Equations \ref{eqn:case1-first-eqn}-\ref{eqn:case1-last-eqn} are satisfied if
\begin{equation}\label{lemma-equation-1}\th_{p,n}(q)=\th_{p,b}(q) \text{ and }  \th_{p,2m+n}(q)=\th_{p,2a+b}(q).
\end{equation}
Equations \ref{eqn:case2-first-eqn}-\ref{eqn:case2-last-eqn} are satisfied if
\begin{equation}\label{lemma-equation-2}\th_{p,n}(q)=\th_{p,b+p}(q) \text{ and } \th_{p,2m+n}(q)=\th_{p,2a+b+p}(q).
\end{equation}
That is, if Equation \ref{lemma-equation-1} or \ref{lemma-equation-2} holds, then $\th_{\La_{m,n}}(q)=\th_{\La_{a,b}}(q)$.

From Eq.~\ref{lemma-equation-1}, we have four subcases corresponding to $n\equiv\pm b\mod 2p$ and
$2m+n\equiv\pm(2a+b)\mod 2p$.  If $n\equiv b\mod 2p$, one finds that $m\equiv a\mod p$ or $m\equiv -a-b \mod
p$. If $n\equiv -b\mod 2p$, one finds that $m\equiv a+b\mod p$ or $m\equiv -a \mod p$.

From Eq.~\ref{lemma-equation-2}, we have four subcases as well, corresponding to $n\equiv \pm(b+p)\mod 2p$
and $2m+n\equiv\pm(2a+b+p)\mod 2p$. If $n\equiv b+p \mod 2p$, then either $m\equiv a\mod p$ or $m\equiv -a-b
\mod p$.  And if $n\equiv -b-p \mod 2p$, then either $m\equiv a+b\mod p$ or $m\equiv -a \mod p$. Therefore,
if $n\equiv b \mod p$, then $m\equiv a \mod p$ or $m\equiv -a-b \mod p$. If $n\equiv -b\mod p$, then $m\equiv
a+b\mod p$ or $m\equiv -a\mod p$. \qed

The Klein 4-group  generated by matrices
\[ \begin{pmatrix}-1&0\\0&-1\end{pmatrix} \text{   and   }\begin{pmatrix}1&1\\0&-1\end{pmatrix}\]
acts on $(Z/pZ)^2$. The orbits form equivalence classes  on $\Z^2$. This equivalence is given by
\[(a,b)\sim(m,n) \text{ if } (m,n) \equiv  (a,b), \, (-a-b, b), \, (-a, -b), \text{ or } (a+b, -b) \mod p.\]
By Lemma \ref{lemma-congruences}, if $(a, b)\sim (m, n)$, then \[\th_{\La_{a,b}}(q)=\th_{\La_{m,n}}(q).\]
Then we have the following result:
\begin{cor}\label{cor-num-theta-fns}
For   any odd prime $p$, the set $\set{\th_{\La_{a,b}}(q) : a, b \in \Z}$ contains at most
$\frac{(p+1)^2}{4}$ elements.
\end{cor}

\proof We will prove this by showing that there are $\frac{(p+1)^2}{4}$ equivalence classes under the
relation $\sim$.  This will imply that there are at most $\frac{(p+1)^2}{4}$ theta functions.  Note that
$(a,b)\sim(a+mp, b+np)$ for any $m, n\in\Z$.  Thus, it is enough to consider only $a, b \in \set{0, \dots,
p-1}$.

Consider the equivalence class of $(a,b)$, which is $$\set{(a,b), (-a-b, b), (-a, -b), (a+b, -b)}.$$  This
set contains either 1, 2, or 4 elements.  (If two elements are equal, then the two remaining elements are
also equal.)  If $b=0$, the set contains $(a, 0)$ and $(-a, 0)$, which are equal if $a=0$ and non-equal if
$a\neq0$ (using the fact that $p$ is odd).  Thus, if $b=0$, there is one equivalence class corresponding to
$a=0$ and there are $\frac{p-1}{2}$ equivalence classes containing elements with $a\neq 0$.

If $b\neq0$, then $b\not\equiv -b \mod p$, so $(a,b)$ and $(-a, -b)$ are distinct mod $p$.  This means there
are either 2 or 4 elements in the equivalence class of $(a,b)$.  Further, $(a,b)$ and $(-a-b,b)$ are
congruent mod $p$ if and only if $(-a, -b)$ and $(a+b, -b)$ are congruent mod $p$ if and only if $2a\equiv
-b\mod p$.  Thus, if $2a\equiv -b\mod p$, the equivalence class of $(a,b)$ has 2 elements.  There are $p-1$
pairs $(a,b)$ with $b\neq0$ and $2a\equiv -b\mod p$, which gives $\frac{p-1}{2}$ equivalence classes.  There
are $(p-1)^2$ remaining pairs $(a,b)$ for which $b\neq0$ and $2a\not\equiv-b\mod p$.  The equivalence classes
for these pairs contain 4 elements, leading to $\frac{(p-1)^2}{4}$ equivalence classes. Summed up, we have
$1+\frac{p-1}{2}+\frac{p-1}{2}+\frac{(p-1)^2}{4}=\frac{(p+1)^2}{4}$ equivalence classes, meaning there are at
most $\frac{(p+1)^2}{4}$ theta functions. \qed

The next result determines in what cases we have exactly $\frac{(p+1)^2}{4}$ theta functions.

\begin{thm}\label{bound_on_d_for_dimension}
For any odd prime $p$ and any  $d>3p^2$, the set $\set{\th_{\La_{a,b}}(q) : a,b\in\Z}$ spans a
$\frac{(p+1)^2}{4}$ dimension vector space in $\Z[[q]]$.  Hence,  Lemma \ref{lemma-congruences} is an ``if and only  if'' statement for large enough $d$.
\end{thm}
\proof We prove this by calculating the minimal exponent appearing in the power series of
$\th_{\La_{a,b}}(q)$ for any $a,b\in\Z$.  We will find that there are $\frac{(p+1)^2}{4}$ different such
minimal exponents, indicating that there is no linear relationship between the $\frac{(p+1)^2}{4}$
corresponding theta series.  From Corollary \ref{cor-num-theta-fns}, there are at most $\frac{(p+1)^2}{4}$ such
series, so we can then conclude that there are exactly $\frac{(p+1)^2}{4}$ such series. Let $a,b\in\Z$ with
$0\leq a<p$ and $0\leq b < p$.  Expanding $\th_{\La_{a,b}}(q)$, one finds that
\begin{equation*}
\begin{split}
\th_{\La_{a,b}}(q)   = & ~ \th_{p,b}(q^{p^2\ell})\th_{p,2a+b}(q^{p^2})+\th_{p,b+p}(q^{p^2\ell})\th_{p,2a+b+p}(q^{p^2}) \\ \\
  = & ~\sum_{n\in\Z} q^{p^2\ell(n+b/2p)^2}\sum_{m\in\Z} q^{p^2(m+a/p+b/2p)^2} \\
 &  + \sum_{n\in\Z} q^{p^2\ell(n+1/2+b/2p)^2}\sum_{m\in\Z} q^{p^2(m+1/2+a/p+b/2p)^2} \\ \\
 = & ~\sum_{n\in\Z} q^{(\ell/4)(2pn+b)^2}\sum_{m\in\Z} q^{(1/4)(2pm+2a+b)^2} \\
&  + \sum_{n\in\Z} q^{(\ell/4)(2pn+p+b)^2}\sum_{m\in\Z} q^{(1/4)(2pm+p+2a+b)^2}.
\end{split}
\end{equation*}
Using the fact that $0\leq b < p$, the term with the smallest exponent in the first summation is
$q^{(\ell/4)b^2}$ and the term with the smallest exponent in the second summation is either
$q^{(1/4)(2a+b)^2}$ or $q^{(1/4)(2a+b-2p)^2}$ (depending on how big $2a+b$ is).  Thus, the term with minimal
exponent in the product of the first two summations is either $$q^{(\ell/4)b^2}\cdot q^{(1/4)(2a+b)^2}\text{
or } q^{(\ell/4)b^2}\cdot q^{(1/4)(2a+b-2p)^2}.$$  Using the fact that $\ell=4d-1$, this term is either
\[q^{a^2+ab+db^2}\text{ or }q^{(a-p)^2+(a-p)b+db^2}.\]
Working analogously with the product of the second pair of summations, one finds the term with smallest
exponent there is either
\[q^{a^2+a(b-p)+d(b-p)^2}\text{ or }q^{(a-p)^2+(a-p)(b-p)+d(b-p)^2}.\]
Thus, in the theta series $\th_{\La_{a,b}}(q)$, the smallest power of $q$ is the minimum of
\[
\begin{split}
& a^2+ab+db^2,  \quad (a-p)^2+(a-p)b+db^2,\\
&  a^2+a(b-p)+d(b-p)^2, \quad  (a-p)^2+(a-p)(b-p)+d(b-p)^2.
\end{split}
\]
Let $min(\th_{\La_{a,b}}(q))$ denote this minimal exponent. Suppose that
$min(\th_{\La_{a,b}}(q))=min(\th_{\La_{m,n}}(q))$ for some integers $a, b, m, n\in\set{0,1,\dots,p-1}$ and
some value of $d$.  Then, $min(\th_{\La_{a,b}}(q))=u^2+uv+dv^2$, where $u=a$ or $u=a-p$ and $v=b$ or $v=b-p$.
Similarly, $min(\th_{\La_{m,n}}(q))=x^2+xy+dy^2$ where $x=m$ or $x=m-p$ and $y=n$ or $y=n-p$.  Note that we
have $|u|, |v|, |x|, |y| \leq p$. We have two cases to consider, either $v^2\neq y^2$ or $v^2=y^2$.

If $v^2\neq y^2$, then, solving for $d$, we find that
$d=\frac{u^2+uv-x^2-xy}{y^2-v^2}.$
Thus,
\[
|d|  \leq  \frac{|u^2-x^2|+|uv|+|-xy|}{|y^2-v^2|} \leq  \frac{p^2+p^2+p^2}{|y^2-v^2|} \leq 3p^2.
\]
If $v^2=y^2$, then given that $u^2+uv+dv^2=x^2+xy+dy^2$, we find that $y=\pm v$.  If $y=v$, then we find
$u^2+uv=x^2+xv$, so $u^2-x^2+uv-xv=0$, so $(u-x)(u+x+v)=0$.  Thus, $x=u$ or $x=-u-v$.  Similarly, if $y=-v$,
then $u^2+uv=x^2-xv$, which implies that $x=-u$ or $x=u+v$.

Using the facts that $u\equiv a \mod p$, $v\equiv b \mod p$, $x\equiv m \mod p$, and $y\equiv n \mod p$, we
find that $(m,n)$ is congruent modulo $p$ to one of the ordered pairs $(a, b)$, $(-a-b, b)$, $(-a, -b)$,
$(a+b, -b)$.
Hence, if $d>3p^2$, then $\th_{\La_{a,b}}(q)=\th_{\La_{m,n}}(q)$ if and only if $(a,b)\sim(m,n)$.  By the
above corollary, there are precisely $\frac{(p+1)^2}{4}$ equivalence classes. Hence, there are precisely
$\frac{(p+1)^2}{4}$ theta functions.  Furthermore, since these theta functions all have different leading
exponents, they are linearly independent. This completes the proof. \qed

\begin{rem}
The bound for $d$ given in Theorem \ref{bound_on_d_for_dimension} is not sharp.  For instance, using a computer algebra package, one finds that for $d=2$, there are $\frac{(p+1)^2}{4}$ equivalence classes for all primes $p\leq19$.
\end{rem}

\section{Theta functions of codes over $\O_K / p\O_K$}
Let $p\nmid\ell$ and $\R := \O_K / p\O_K =\set{a+b\om : a,b\in\F_p, \om^2+\om+d=0}.$ We have the map
\[\rho_{\ell, p}:\O_K \rightarrow \O_K/ p\O_k =: \R\]
A linear code $C$ of length $n$ over $\R$ is an $\R$-submodule of $\R^n$. The dual is defined as
$C^\bot=\set{u\in \R^n : u\cdot \bar{v}=0 \text{ for all } v \in C}$. If $C=C^\bot$ then $C$ is self-dual.
We define
\[\L_{\ell}(C):= \set { u=(u_1,\dots,u_n) \in \O_K^n : (\rho_{\ell, p}(u_1),\dots,\rho_{\ell,p}(u_n)) \in C},\]
In other words, $\L_{\ell}(C)$ consists of all vectors in $\O_K^n$ in the inverse image of $C$, taken componentwise by $\rho_{\ell,p}$.
This method of lattice construction is known as Construction A.

For notation, let $r_{a+pb+1}=a-b\om$, so $\R=\set{r_1,\dots,r_{p^2}}.$  For a codeword
$u=(u_1,\dots,u_n)\in\mathcal{R}^n$ and $r_i \in\mathcal{R}$, we define the counting function $n_i(u)
:=\#\set{i : u_i=r_i}.$  The complete weight enumerator of the $\R$ code $C$ is the polynomial
\begin{equation}
cwe_C(z_1, z_2, \dots, z_{p^2})=\sum_{u\in C}z_1^{n_1(u)}z_2^{n_2(u)}\dots z_{p^2}^{n_{p^2}(u)}.
\end{equation}
We can use this polynomial to find the theta function of the lattice $\La_\ell(C)$. For a proof of the
following result see \cite{SS}.
\begin{lem}\label{lem1} Let $C$ be a code defined over $\R$ and $cwe_C$ its complete weight enumerator as above. Then,
\[\th_{\L_\ell(\mathcal{C})}(q) = cwe_\mathcal{C}(\th_{\La_{0,0}}(q),\th_{\La_{1,0}}(q),\dots,\th_{\La_{p-1,p-1}}(q))\]
\end{lem}

In \cite{ch}, for $p=2$, the symmetric weight enumerator polynomial $swe_\mathcal{C}$ of a code $\mathcal{C}$
over a ring or field of cardinality 4 is defined to be $$swe_\mathcal{C}(X,Y,Z)=cwe_\mathcal{C}(X,Y,Z,Z).$$
For $\La_{\Lambda_\ell ( \mathcal{C} )}(q)$, the lattice obtained from $\mathcal{C}$ by Construction A, by Theorem 5.2 of
\cite{ch}, one can then write
$$\th_{\La_\ell(\mathcal{C})}(q)=swe_\mathcal{C}(\th_{\La_{0,0}}(q),\th_{\La_{1,0}}(q),\th_{\La_{0,1}}(q)).$$
These theta functions are referred to as $A_d(q), C_d(q),$ and $G_d(q)$ in \cite{ch} and \cite{SV}.

\begin{rem}
The connection between complete weight enumerators of self-dual codes  over $\F_p$ and Siegel theta series of
unimodular lattices is well known. Construction A associates to any length $n$ code $C=C^\perp$ an
$n$-dimensional unimodular lattice; see \cite{LS} for details.

\end{rem}

For $p>2$, however, there are $\frac{(p+1)^2}{4}$   theta functions associated to the various lattices, so
our analog of the symmetric weight enumerator polynomial needs more than 3 variables.

\begin{exmp}
For $p=3$, from Remark 2.2 in \cite{SS}, we have four theta functions corresponding to the lattices
$\La_{a,b}$, namely $\th_{\La_{0,0}}(q)$, $\th_{\La_{1,0}}(q)$, $\th_{\La_{1,1}}(q)$, $\th_{\La_{0,1}}(q).$
If we define the symmetric weight enumerator to be
$$swe_\mathcal{C}(X,Y,Z,W)=cwe_\mathcal{C}(X,Y,Y,Z,W,Z,Z,Z,W),$$ then by replacing
\[X=\th_{\La_{0,0}}(q), \, \, Y=\th_{\La_{1,0}}(q),  \, \, Z=\th_{\La_{1,1}}(q), \, \, W=\th_{\La_{0,1}}(q),\]
one finds that
\[
\begin{split}
\th_{\La_\ell (C) } (q) \,  =  \, & cwe_\mathcal{C}(\th_{\La_{0,0}}(q),\th_{\La_{1,0}}(q),\dots,\th_{\La_{2,2}}(q)), \\
 \, =\,  & swe_\mathcal{C}(\th_{\La_{0,0}}(q),\th_{\La_{1,0}}(q),\th_{\La_{1,1}}(q),\th_{\La_{0,1}}(q)).
\end{split}
\]
\end{exmp}

\begin{exmp}
Let $\mathcal{C}_1$ be the length-2 repetition code $\set{(x,x):x\in\R}$ for $p=3$.  The complete weight
enumerator of this code is $$cwe_{\mathcal{C}_1}(z_0,\dots,z_8)=z_0^2+z_1^2+\dots+z_8^2,$$ and the symmetric
weight enumerator is $$swe_{\mathcal{C}_1}(X,Y,Z,W)=X^2+2Y^2+4Z^2+2W^2.$$  Using some computational algebra
package, one then finds $\th_{\La_\ell{(C)}}(q)$ for each $\ell$. We display the cases when $\ell=7, 11$.
\[
\begin{split}
\th_{\La_7 (C)  }(q) =  &   1 + 2 q^2 +4q^4 + 4 q^5+10 q^8 +4 q^9+16 q^{10}+ 8q^{11}+8q^{13}+2q^{14}\\
& + 24 q^{16}+12q^{17} +12 q^{18}+ 16 q^{19} +28 q^{20}+20q^{22} +16q^{23} +16q^{25}\\
& + 28q^{26}+ 16q ^{27} +4 q^{28} + 20 q^{29} +24 q^{31}+ 42q^{32}+ 32q^{34} +4q^{35} \\
& + 28 q^{36}+ 24q^{37} + 40 q^{38} +56q^{40}+28 q^{41} + 32q^{43}+56q^{44}+ 24q^{45}\\
& + 52q^{46}+ 32q^{47} + 62q^{50}+ \cdots\\
\th_{\La_{11}(C) }(q)  = & 1+2\,{q}^{2}+4\,{q}^{5}+4\,{q}^{6}+2\,{q}^{8}+4\,{q
}^{9}+8\,{q}^{10}+8\,{q}^{12}+8\,{q}^{15}
+8\,{q}^{16}\\
&+4\,{q}^{17} +24\,{q}^ {18}+4\,{q}^{20}+8\,{q}^{21}+2\,{q}^{22}+20\,{q}^{24}+ 16\,{q}^{25}+12\,{q}^{26}\\
&+24\,{q}^{27}+16\,{q}^{28}+12\,{q}^{29}+24\,{q}^{30}+8\,{q}^{
31}+10\,{q}^{32}+8\,{q}^{34}+8\,{q}^{35}\\
&+36\,{q}^{36}+8\,{q}^{38}+16\,{q}^{39}+8\,{q}^{40}+20\,{q}^{41}+24\,{q}^{42}+16\,{q}^{43}+32\,{q}^ {45}\\
&+8\,{q}^{46}+8\,{q}^{47}+40\,{q}^{48}+8\,{q}^{49}+22\,{q}^{50}+\cdots\\
\end{split}
\]
It will be the goal of our next section to study how the corresponding theta function of a given code differ
for different levels $\ell$.
\end{exmp}

In general, finding such an explicit relation between the theta function and the symmetric weight enumerator
polynomial for larger $p$ seems difficult.
\begin{prob}
Determine an explicit relation between theta functions and the symmetric weight enumerator polynomial of a
code defined over $\R$ for $p> 3$.
\end{prob}
We expect that the answer to the above problem is that the theta function is given as the symmetric weight
enumerator $swe_C $ of $C$, evaluated on the  theta functions defined on cosets of $\O_K/p\O_K$.

\section{Theta functions and the corresponding complete weight enumerator polynomials}
For a fixed prime $p$, let $C$ be a linear code over $\R=\F_{p^2}$ or $\F_p\times\F_p$ of length $n$ and
dimension $k$. An admissible level $\ell$ is an integer $\ell$ such that $\O_K/p\O_K$ is isomorphic to
$\R$.  For an admissible $\ell$, let $\L_{\ell} (C) $ be the corresponding lattice as in the previous section.
Then, the \textbf{level $\ell$ theta function} $\Th_{\L_{\ell} (C) } (\tau) $ of the lattice $\L_{\ell} (C)$
is determined by the complete weight enumerator $cwe_C $ of $C$, evaluated on the  theta functions defined on
cosets of $\O_K/p\O_K$.  We consider the following questions. How do the theta functions $\Th_{\L_{\ell} (C)
} (\tau) $ of the same code $C$ differ for different levels $\ell$? Can non-equivalent codes give  the same
theta functions for all levels $\ell$?

We give a satisfactory answer to the first question (cf. Theorem~\ref{thm1}, Lemma~\ref{lem2})  and for the
second question we conjecture that:
\begin{conj}\label{conj1}
Let   $C$ be a code of size $n$ defined over $\R$ and $\Th_{\La_\ell (C)}$ be its corresponding theta
function for level $\ell$. Then, for large enough $\ell$, there is a unique complete weight enumerator
polynomial which corresponds to $\Th_{\La_\ell (C)}$.
\end{conj}

Let $C$ be a code defined over $\R$ for a fixed $p >2$. Let the complete weight enumerator of $C$ be the
degree $n$ polynomial $cwe_C = f( x_1, \dots , x_r)$, for $r= p^2$. Then from Lemma~\ref{lem1} we have that
\[\th_{\L_\ell(\mathcal{C})}(\tau) = f(\th_{\La_{0,0}}(\tau),\dots,\th_{\La_{p-1,p-1}}(\tau))\]
for a given $\ell$. First we want to address how $\th_{\L_\ell(\mathcal{C})}(\tau)$ and
$\th_{\L_{\ell^\prime} (\mathcal{C})}(\tau)$ differ for different $\ell$ and $\ell^\prime$. The proof of the
following remark is elementary.

\begin{rem}
For $n\neq0$, $Q_d(m,n)\geq d$.
\end{rem}

Then we have  the following theorem.

\begin{thm}\label{thm1}
Let $C$ be a code defined over $\R$.  For all admissible $\ell, \ell'$ with $\ell<\ell'$ the following holds
\[\th_{\La_\ell(C)}(q)=\th_{\La_{\ell'}(C)}(q)+\O(q^{\frac {\ell+1} 4}).\]
\end{thm}

\proof From section 3, we have the map $\rho_{\ell,p}:\O_K\rightarrow\O_K/p\O_K\rightarrow\R$ and
\[\La_\ell(C)=\set{u=(u_1, \dots, u_n)\in\O_K^n : (\rho_{\ell,p}(u_1), \dots, \rho_{\ell,p}(u_n))\in C}.\]
We denote $u_i=a_i-b_i\om_\ell$ for $a_i, b_i \in\Z$ with $i=1, \dots, n$ and $d=\frac {\ell+1} 4$. Then
\begin{equation*}
\begin{split}
\th_{\La_\ell(C)}(q) & =
\displaystyle\sum_{u\in\La_\ell(C)}q^{u\cdot \bar{u}}  =  \displaystyle\sum_{u\in\La_\ell(C)}q^{u_1\bar{u_1}+\dots+u_n\bar{u_n}} \\ \\
& =  \displaystyle\sum_{u\in\La_\ell(C)}q^{Q_d(a_1, b_1)+\dots+Q_d(a_n, b_n)} \\ \\
& = \displaystyle\sum_{\substack {u\in\La_\ell(C), \\ b_i=0 \\ \text{ for all } i} }q^{Q_d(a_1,
b_1)+\dots+Q_d(a_n, b_n)}
+ \displaystyle\sum_{\substack {u\in\La_\ell(C), \\ b_i\neq0 \\ \text{ for some } i}} q^{Q_d(a_1, b_1)+\dots+Q_d(a_n, b_n)} \\ \\
& =  \displaystyle\sum_{\substack  {u\in\La_\ell(C), \\ b_i=0 \text{ for all } i}  }q^{a_1^2+\dots+a_n^2} +
\displaystyle\sum_{\substack{ u\in\La_\ell(C),\\ b_i\neq0 \text{ for some } i } }q^{Q_d(a_1,
b_1)+\dots+Q_d(a_n, b_n)}.
\end{split}
\end{equation*}

\noindent  Note that this first summation does not depend on $d$ (or $\ell$). In the second summation, each
term's exponent contains a term of the form $Q_d(a_i, b_i)$ where $b_i\neq0$. By the lemma above, we
have $Q_d(a_i, b_i)\geq d$.  Since all of the terms in the exponent are added, each term in the second
summation has exponent at least $d$.  Hence, the second summation is $\O(q^d)$. Thus, we have
\[\th_{\La_\ell(C)}(q)=\displaystyle\sum_{\substack{u\in\La_\ell(C), \\ b_i=0 \text{ for all } i}}q^{a_1^2+\dots+a_n^2}+\O(q^d).\]
Similarly,
\[\th_{\La_{\ell'}(C)}(q)=\displaystyle\sum_{\substack{u\in\La_{\ell'}(C),\\ b_i=0 \text{ for all }i}}q^{a_1^2+\dots+a_n^2}+\O(q^{d'}).\]
For admissible $\ell$, $\ell'$ with $\ell<\ell'$, we conclude that
$$\th_{\La_\ell(C)}(q)-\th_{\La_{\ell'}(C)}(q)=\O(q^d).$$
This completes the proof. \qed

We have the following lemma.
\begin{lem}\label{lem2}
Let $C$ be a fixed code of size $n$ defined over $\R$ and $\Th(q)=\sum \lambda_i q^i$ be its theta
function for level $\ell$. Then, there exists a bound $B_{\ell,p,n}$ such that  $\Th(q)$ is uniquely
determined by its first $B_{\ell,p,n}$ coefficients.
\end{lem}

\proof  We want to show that if $\Th(q)-\Th'(q)=\O(q^{B_{\ell,p,n}})$, then $\Th(q)=\Th'(q)$.  Fix $p, n,
\ell$. There are finitely many codes $C$ over $\R$ of length $n$.  Denote them by $C_1, \dots, C_m$, for
some integer $m$.  To each code $C_i$, there is a corresponding theta function $\Th_{C_i}(q)$.

Let $$S=\set{r\in\Z_{\geq0} : \Th_{C_i}(q)-\Th_{C_j}(q)=\O(q^r) \text{ and }
\Th_{C_i}(q)\neq\Th_{C_j}(q)}$$ and let $B_{\ell,p,n}=1+\max{S}$.  Since $S$ is finite, $B_{\ell,p,n}$
is well-defined. Furthermore, if $\Th_{C_i}(q)-\Th_{C_j}(q)=\O(q^m)$ for some $m\geq B_{\ell,p,n}$,
this implies that $m\notin S$, so we must have $\Th_{C_i}(q)=\Th_{C_j}(q)$. \qed

For notation, when $p$ and $n$ are fixed, we will let $B_\ell=B_{\ell,p,n}$.

To extend the theory for $p=2$ to $p>2$ we have to find a relation between the theta function $\Th_{\La_\ell
(C)}$ and the number of complete weight enumerator polynomials corresponding to it.

Fix an odd prime $p$ and let $C$ be a given code of length $n$ over $\R$. Choose an admissible value of
$\ell$ such that there are $\frac{(p+1)^2}{4}$ independent theta functions (as in Theorem
\ref{bound_on_d_for_dimension}).  Then, the complete weight enumerator of $C$ has degree $n$ and
$r=\frac {(p+1)^2} 4$ variables $x_1, \dots , x_r$. We call a \emph{generic complete weight enumerator
polynomial} a homogenous polynomial  $P \in \Q [x_1, \dots , x_r]$.

\begin{lem}
A degree $n$ generic complete weight enumerator polynomial  has   $s:=\frac {\left( n-1 + \frac {(p+1)^2 } 4
\right) \text{!}} {{n \text!} \cdot  \left(  \frac {(p+1)^2} 4 -1 \right) \text{!} }$ monomials.
\end{lem}

\proof We need to count the number of monomials of a homogenous degree $n$ polynomial in $r=\frac {(p+1)^2 }
4$ variables which is
\[ s=\frac {(n+r-1) \text{!} } {n \text{!}  \, (r-1) \text{!} } =
\frac {\left( n-1 + \frac {(p+1)^2 } 4 \right) \text{!}} {{n \text!} \cdot  \left(  \frac {(p+1)^2} 4 -1
\right) \text{!} }
\]
This completes the proof. \qed

Denote by   $P(x_1, \dots , x_r )$  a generic $r$-nary, degree $n$, homogeneous polynomial.  Assume that
there is a length $n$ code $C$ defined over $\R$ such that $P(x_1, \dots , x_r )$ is the symmetric weight
enumerator polynomial. In other words,
\[swe_C (x_1,\dots, x_r)=P(x_1, \dots , x_r)\]
Fix the level $\ell$. Then, by replacing
\[x_1 = \Th_{\Lambda_{0,0}}(q),   \, \dots \dots \,  ,   x_r=\Th_{\Lambda_{p-1,p-1}}(q),\]
we compute the left side of the above equation as a series $\sum_{i=0}^\infty \lambda_i q^i$. By equating
both sides of $ \sum_{i=0}^\infty \lambda_i q^i = P(x_1, \dots , x_r ), $
we can get a linear system of equations.  Since the first $\lambda_0, \dots, \lambda_{B_\ell-1}$ determine all
the coefficients of the theta series then we have to pick $B_\ell$ equations (these equations are not
necessarily independent).

Consider the coefficients of the polynomial $P(x_1, \dots , x_r)$ as parameters $c_1, \dots c_s$. Then, the
linear map
\[
\begin{split}
 L_{\ell}:  \C^s  &  \to \C^{B_\ell-1}  \\
  (c_1, \dots c_s ) & \mapsto (\lambda_0, \dots , \lambda_{B_\ell-1})
 \end{split}
 \]
has an associated matrix $M_{\ell}$.  For a fixed value of  $(\lambda_0, \dots , \lambda_{B_\ell-1})$,
determining the rank of the matrix $M_{\ell}$ would determine the number of polynomials giving the same
theta series. There is a unique complete weight enumerator corresponding to a given theta function when
\[ \mbox{null } ( M_{\ell}) = s - \mbox{rank } ( M_{\ell} ) =0\]

\begin{conj}\label{conj2} For $\ell \geq \frac {p(n+1)(n+2)}{n} -1$ we have $\mbox{\emph{null} } M_{\ell} =0$, or in other
words
\[ \mbox{\emph{rank} } ( M_{\ell} )= \frac {\left( n-1 + \frac {(p+1)^2 } 4 \right) \text{!}} {{n \text!} \cdot  \left(  \frac {(p+1)^2} 4 -1
\right) \text{!} }\]
\end{conj}

The choice of $\ell $ is taken from experimental results for primes $p=2$ and 3. More details are given in
the next section.

It is obvious that Conjecture~\ref{conj2} implies Conjecture~\ref{conj1}. If Conjecture~\ref{conj1} is true
then for large enough $\ell$ there would be a one to one correspondence between the complete weight
enumerator polynomials and the corresponding theta functions.  Perhaps, more interesting is to find $\ell$ and
$p$ for which there is not a one to one such correspondence.  Consider the map
\[ \Phi ( \ell, p ) = \left( \lambda_0 (\ell, p) , \dots , \lambda_{B_\ell-1} (\ell, p) \right), \]
where $\lambda_0, \dots , \lambda_{B_\ell-1}$ are now functions in $\ell$ and $p$. Let $V$ be the variety given
by the Jacobian of the map $\Phi$. Finding integer points $\ell, p$ on this variety such that $\ell$ and $p$
satisfy our assumptions would give us values for $\ell, p$ when the above correspondence is not one to one.
However, it seems quite hard to get explicit description of the map $\Phi$. Next, we will try to shed some
light over the above  conjectures for fixed small primes $p$.

\section{Bounds for small primes}
In \cite{SV} we determine explicit bounds for the above theorems for prime $p=2$. In this section we give some
computation evidence for the generalization of the result for $p=3\ .$ We recall the theorem for $p=2.$
\begin{thm}[\cite{SV}, Thm.~2]
Let $p=2$ and $C$ be a code of size $n$ defined over $\R$ and $\Th_{\La_\ell} (C)$ be its corresponding theta
function for level $\ell$. Then the following hold:
\begin{description}
\item [i)]
For $\ell < \frac {2(n+1)(n+2)}{n} -1$ there is a $\delta$-dimensional family of symmetrized weight
enumerator polynomials corresponding to $\Th_{\La_\ell}(C)$, where\\ $\delta \geq \frac
{(n+1)(n+2)}{2}-\frac{n(\ell+1)}{4} - 1$.

\item [ii)]
For $\ell \geq \frac {2(n+1)(n+2)}{n} -1$ and $n < \frac{\ell+1}{4}$ there is a unique symmetrized weight
enumerator polynomial which corresponds to $\Th_{\La_\ell}(C)$.
\end{description}
\end{thm}

These results were obtained by using the explicit expression of theta in terms of the symmetric weight
enumerator valuated on the theta functions of the cosets.

Next we want to find explicit bounds for $p=3$ as in the case of $p=2$. In the case of $p=3$ it is enough to consider
four theta functions, $\th_{\Lambda_{0,0}}(q)$, $\th_{\Lambda_{1,0}}(q)$, $\th_{\Lambda_{0,1}}(q),$ and
$\th_{\Lambda_{1,1}}(q)$ since $\th_{\Lambda_{2,0}}(q) = \th_{\Lambda_{1,0}}(q)$\,\,,\,\,$\th_{\Lambda_{2,2}}(q) =
\th_{\Lambda_{1,1}}(q)$ and $\th_{\Lambda_{0,2}}(q) = \th_{\Lambda_{1,2}}(q) = \th_{\Lambda_{2,1}}(q) =
\th_{\Lambda_{0,1}}(q).$ If we are given a code $C$ and its weight enumerator polynomial then we can find the theta
function of the lattice constructed from $C$ using Construction A. Let $\Th(q) = \sum_{i=0}^\infty \lambda_i q^i$
be the theta series for level $\ell$ and
\[p(x,y,z,w)=\sum_{i+j+k+m=n}    c_{i,j,k}x^iy^jz^k w^m\]
be a degree $n$ generic 4-nary homogeneous polynomial. We would like to find out how many polynomials $p(x,y,z,w)$
correspond to $\Th(q)$ for a fixed $\ell$. For a given $\ell$ find $\th_{\Lambda_{0,0}}(q),$ $\th_{\Lambda_{1,0}}(q),$
$\th_{\Lambda_{0,1}}(q)$ and $\th_{\Lambda_{1,1}}(q)$ and substitute them in the $p(x,y,z,w).$ Hence, $p(x,y,z,w)$ is
now written as a series in $q.$ We get infinitely many equations by equating the corresponding coefficients of the two
sides of the equation
\[ p(\th_{\Lambda_{0,0}}(q), \th_{\Lambda_{1,0}}(q),
\th_{\Lambda_{0,1}}(q), \th_{\Lambda_{1,1}}(q) ) = \sum_{i=0}^\infty \lambda_i q^i.\]
Since the first $\lambda_0, \dots, \lambda_{B_\ell-1}$ determine all the coefficients of the theta series then it is
enough to pick the first $B_\ell$ equations. The linear map
\[ L_{\ell}: (c_1, \dots c_{20} ) \mapsto (\lambda_0, \dots , \lambda_{B_\ell-1})\]
has an associated matrix $M_{\ell}$. If the nullity of $M_{\ell}$ is zero then we have a unique polynomial
that corresponds to the given  theta series. We have calculated the nullity of the matrix and $B_\ell$ for
small $n$ and $\ell.$

\begin{exmp}[The case $p=3, n=3$]
The generic homogenous polynomial is given by
\begin{equation}
\begin{split}
P(x,y,z)& =c_1 x^3 + c_2 x^2y + c_3 x^2z + c_4x^2w + c_5xy^2 + c_6xz^2 + c_7xw^2 + c_8xyz \\
       & + c_9xyw + c_{10}xzw + c_{11}y^3 + c_{12}y^2z + c_{13}y^2w + c_{14}yz^2 + c_{15}yw^2 \\
       & + c_{16}yzw + c_{17}z^3 + c_{18}z^2w + c_{19}zw^2 + c_{20}w^3.\\
\end{split}
\end{equation}
The system of equations can be written by the form of \[A \vec{c} = \vec{\lambda}\] where $\vec{c} =
\begin{pmatrix} c_1 &  c_2 & \cdots & c_{20} \end{pmatrix} ^t$, $\vec{\lambda} =
\begin{pmatrix}\lambda_0 & \lambda_1 & \cdots & \lambda_{15} \end{pmatrix} ^t.$
In the case of $\ell=7$   the matrix $M_7$  has $\mbox{null } ( M_7 ) =4$. We have a positive dimension
family of solution set. The case of $\ell = 11$ the matrix $M_{11}$  has $\mbox{null } (M_{11}) =1$. For any
case where $\ell \geq 19$ the nullity of the matrix is 0. Hence, for every given theta series, there is a
unique symmetric weight enumerator polynomial. .
\end{exmp}

We summarize the results in the following table:\\

\begin{center}
\begin{tabular}{c||c|c||c|c||c|c}
%
\textbf{$\ell$} & \multicolumn{2}{c||}{$n=3$} & \multicolumn{2}{c||}{$n=4$} & \multicolumn{2}{c}{$n=5$}\\
\cline{2-3} \cline{4-5} \cline{6-7}
                & $B_{\ell}$ & $\mbox{\emph{null} } M_{\ell}$ & $B_{\ell}$ & $\mbox{\emph{null} } M_{\ell}$
                & $B_{\ell}$ & $\mbox{\emph{null }} M_{\ell}$ \\
 \hline
7 & 16 & 4  & 26 & 9 & 33 & 24 \\
\hline
11 & 19 & 1 & 30 & 5  & 42 & 14 \\
\hline
19 & 22 & 0 & 38 & 0 & 60 & 0 \\
\hline
23 & 25 & 0 & 37 & 0 & 58 & 0 \\
\hline
31 & 31 & 0  & 41 & 0  & 60 & 0 \\
\hline
35 & 34 & 0 & 48 & 0 & 61 & 0 \\
\hline
43 & 40 & 0 & 55 & 0  & 69 & 0 \\
\hline
47 & 43 & 0 & 60 & 0 & 74 & 0 \\
\hline
55 & 49 & 0 & 70 & 0 & 86 & 0 \\
\hline
59 & 52 & 0 & 75 & 0 & 92 & 0 \\
%
\end{tabular}
\end{center}

\medskip

Recall that $\ell \equiv 3 \mod 4$ and $p \nmid  \ell$.    It seems  from the table that the same bound of
$B_\ell = \frac {2 (n+1)(n+2)} n$ as for $p=2$ holds also for $p=3, n=3$.

We have the following conjecture for general $p,n$ and $\ell.$

\begin{conj} For a given theta function $\Th_{\La_\ell (C) }$ of a code $C$ for level $\ell$ there is a unique
complete weight enumerator polynomial corresponding to $\Th_{\La_\ell  (C)}$ if $\ell \geq
\frac{p(n+1)(n+2)}{n}.$
\end{conj}

It is interesting to consider such question for such lattices  independently of the connection to coding
theory. What is the meaning of the bound $B_\ell$ for the ring $\O_K/p\O_K$? Do the theta functions
defined here correspond to any modular forms? Is there any difference between the cases when the ring is
$\F_p \times \F_p$ or $\F_{p^2}$?

\noindent{\textbf{Acknowledgments:}}  The authors want to thank the anonymous  referee for useful suggestions.

\end{document}